\documentclass{article}

 \usepackage{amssymb,amsmath}

\begin{document}

\def\Q{{\mathbb Q}}
 \def\R{{\mathbb R}}
 \def\C{{\mathbb C}}
\def\O{{\mathbb O}}
\def\L{{\rm Lat}}
\def\Z{{\mathbb Z}}
\def\F{{\mathbb F}}
\def\Q{{\mathbb Q}}

\def\cV{{\cal V}}
\def\cH{{\cal H}}
\def\cA{{\cal A}}

\def\fP{{\frak P}}
\def\fI{{\frak I}}
\def\fN{{\frak N}}
\def\fF{{\frak F}}

\def\kk{{\bold k}}

\def\phi{\varphi}
\def\kappa{\varkappa}
\def\epsilon{\varepsilon}
\def\le{\leqslant}
\def\ge{\geqslant}

\def\GL{{\rm GL}}
\def\Sp{{\rm Sp}}

\def\Fl{{\rm Fl}}

\def\wt{\widetilde}

\def\tto{\rightrightarrows}

\newcommand{\codim}{\mathop{\rm codim}\nolimits}
\newcommand{\vol}{\mathop{\rm vol}\nolimits}
\newcommand{\res}{\mathop{\rm res}\nolimits}
\renewcommand{\Re}{\mathop{\rm Re}\nolimits}

\newcommand{\Dom}{\mathop{\rm Dom}\nolimits}
\newcommand{\Ker}{\mathop{\rm Ker}\nolimits}
\renewcommand{\Im}{\mathop{\rm Im}\nolimits}
\newcommand{\Indef}{\mathop{\rm Indef}\nolimits}

\newcounter{sec}
\renewcommand{\theequation}{\arabic{sec}.\arabic{equation}}

\begin{center}

{\large\bf

On compression of Bruhat--Tits buildings

}

\medskip

\sc\large Yurii A. Neretin

\end{center}

{\small
Consider
 an affine Bruhat-Tits building $\L_n$
of the type $A_{n-1}$
and the complex distance in $\L_n$, i.e.,
the complete system of invariants of a pair
of vertices of the building.
An elements of the Nazarov semigroup
 is a  lattice
in the duplicated $p$-adic space
 $\Q_p^n\oplus \Q_p^n$.
We investigate behavior
 of the complex distance with respect to
the natural action of the Nazarov semigroup
on the building.
       }

\medskip

It is well known that affine Bruhat--Tits buildings
(see, for instance, \cite{Bro}, \cite{Gar}, \cite{Mac})
are right $p$-adic analogs of Riemannian noncompact
symmetric spaces.  This analogy has no
a priory explanations, but it  exists
on the levels of  geometry,  harmonic analysis,
and special functions.
Some known examples of this parallel are

--- Gindikin--Karpelevich's and
Macdonald's Plancherel formulae, see \cite{Mac};

--- existing of a general theory including
  Hall--Littlewood polynomials
and zonal spherical functions, see \cite{Mac2};

--- similarity in the geometries at infinity
(Semple--Satake type boundaries, Martin boundary etc.),
see \cite{Gui}, \cite{Ner-hinges};

--- parallel between Horn--Klyachko
(see \cite{Ful}) triangle inequality
for eigenvalues of Hermitian matrices%
\footnote{Interpretation of this inequality in terms of geometry
of symmetric spaces is discussed in \cite{Ner-jordan}}
 and the triangle
inequality for complex distance in buildings%
\footnote{ Let $R$, $S$, $T\in\L_n$ (see notations below). To
describe possible collections of $3n$ integers $\{k_j(R,S)\}$,
$\{k_j(R,T)\}$, $\{k_j(S,T)\}$. It can be easily reduced to the
following (unexpectedly nontrivial) question solved in 1960-s by
Ph.Hall and T.Klein. Let $A$, $B$ be finite Abelian $p$-groups.
For which Abelian $p$-groups $C$ there exists a subgroup $A'$
isomorphic to $A$ such that $C/A'\simeq B$, for a solution, see
\cite{Mac2}, II.3.}

--- Existence of a precise parallels of the
harmonic analysis of
the Berezin kernels  and of the matrix beta-functions
 on the  level of buildings%
\footnote{The matrix gamma-function corresponds to
the zeta-function of Tamagawa,
see  \cite{Mac2}, V.4.},
see \cite{Ner-tits}.

\smallskip

In this note, we continue this parallel and obtain
a (quite simple) analog of the compression of angles theorem
(this theorem appeared in \cite{Ner-boson}
as a lemma for boundedness of some integral operators in
 Fock spaces, see its proof in \cite{Ner-book}, 6.3,
and some geometrical discussion in \cite{Ner-krein},
see also \cite{Kou} on its extension to exceptional groups)
for buildings of the type
$A_n$.

\medskip

{\bf 1. Formulation of result}

\smallskip

{\bf 1.1. Notation.}
We denote by $\Q_p$ the $p$-adic field.
By $|\cdot|$, we denote
the norm on $\Q_p$.
Let $\O_p$ be the ring of $p$-adic integers.
Recall that $\O_p$ consists of elements with norm $\le 1$.
By $\Q_p^n$ we denote the coordinate $n$-dimensional
linear space over $\Q_p$.

By $\GL_n(\Q_p)$  we denote the group of invertible
$n\times n$ matrices
over $\Q_p$. By $\GL_n(\O_p)$ we denote the group of matrices
$g\in\GL_n(\Q_p)$ such that all the matrix elements of $g$
and $g^{-1}$ are integer.

\smallskip

{\bf 1.2. Space of lattices.}
A {\it lattice} in $\Q_p^n$ is a compact open $\O_p$-submodule.
 Recall, that any lattice $R$ can be represented
in the form
\begin{equation}
R=\O f_1\oplus \dots\oplus\O f_n
,
\end{equation}
where $f_1$, \dots, $f_n$ is some basis in $\Q_p^n$, see
\cite{Wei-number}, \S II.2.
We denote  by $\L_n$ the space of lattices
in $\Q_p^n$ (this space is countable, and a natural topology
on $\L_n$ is discrete).

The group $\GL_n(\Q_p)$ acts transitively on $\L_n$.
The stabilizer of the lattice
$\O^n_p$
is $\GL_n(\O)$. Thus,
$$\L_n=\GL_n(\Q_p)/\GL_n(\O_p).$$

\smallskip

{\bf 1.3. Complex distance.}

\smallskip

{\sc Theorem} (see, for instance, \cite{Wei-number}, Theorem 2.2.2).
{\it Let $R$, $S$ be two lattices in $\Q_p^n$.

\smallskip

a) There exists a basis $f_1$, \dots, $f_n\in \Q_p^n$
such that
\begin{equation}
R=\sum_{j=1}^n \O_p f_j,\qquad
  S =\sum_{j=1}^n p^{-k_j(R,S)} \O_p f_j
\end{equation}
Integer numbers
$$
k_1(R,S)\ge k_2(R,S)\ge\dots\ge k_n(R,S)
$$
are uniquely determined by the lattices $R$, $S$.

\smallskip

b) If  $S$, $R$ and $S'$, $R'$ be elements of  $\L_n$
satisfying
$$
k_j(R,S)= k_j(R',S')
$$
then there exists $g\in\GL_n(\Q_p)$ such that}
$$
R=g\cdot R',\qquad S=g\cdot S'
$$

We say that $k_j(R,S)$ is the {\it complex distance} between
$S$, $R$
(this structure is the analog of the complex distance or
stationary angles in symmetric spaces, see, for instance,
\cite{Ner-jordan}).

\smallskip

{\bf 1.4. Product of relations.}
Let $V$, $W$ be sets. A relation $H:V\tto W$
is a subset in $V\times W$.

Let $H:V\tto W$, $G:W\tto Y$ be relations.
Their product $G\cdot H:V\tto Y$ is a relation
consisting of all pairs $v\times y\in V\times Y$
for which
there exists $w\in W$ satisfying
$v\times w\in H$, $w\times y\in G$.

\smallskip

{\bf 1.5. Nazarov semigroup.}
An element of the Nazarov semigroup $\Gamma_n$
is a lattice in $\Q_p^n\times \Q_p^n$.
The product in the Nazarov semigroup is the product of
relations.

\smallskip

{\sc Remark.} It is possible to consider elements of
the group $\GL_n(\Q_p)$ as points of the Nazarov
semigroup  at infinity. Indeed, for $g\in\GL_n(\Q_p)$
consider a sequence of lattices
$Z_j\subset \Q_p^n\oplus\Q_p^n$ such that
$$
\cap_{j=1}^\infty
\cup_{k=j}^\infty
Z_j
$$
coincides with the graph of $g$.
For each lattice $R\subset\Q_p^n$,
we have $Z_j\cdot R=g\cdot R$
starting from some place $j$.

\smallskip

{\sc Remark.} Nazarov semigroups
are $p$-adic analogues of semigroups
of linear relations (see \cite{Ner-book}).
 The main nontrivial property of the Nazarov semigroups
is the following statement. The Weil representation of
$\Sp(2n,\Q_p)$ (see \cite{Wei-groups}) admits a canonical
extension to representation of a Nazarov semigroup related to
buildings of the type $C_n$, see \cite{NNO}, \cite{Naz}. It seems
that this statement must have many corollaries (due to
possibilities of the Howe duality), these corollaries never were
discussed in literature.

\smallskip

{\bf 1.6. Action of the Nazarov semigroup on $\L_n$.}
Consider $H\in\Gamma_n$, $R\in \L_n$, We define $HR\in \L_n$
as set of all $w\in \Q_p^n$ such that there exists
$v\in R$ satisfying $v\times w\in H$.

\smallskip

{\bf 1.7. Result of the paper.}
The main result of this paper is the following  theorem.

\smallskip

{\sc Theorem.}
{\it Let $H$ be an element of the Nazarov semigroup $\Gamma_n$.
Let $R$, $S\in \L_n$. Then for each $j$ the number
$k_j(H\cdot R,H\cdot S)$ lies between $0$ and
$k_j( R, S)$, i.e.,}
\begin{align*}
k_j(R,S)\ge 0
\qquad \text{implies}\qquad
 0\le k_j(H\cdot R,H\cdot S)\le k_j(R,S) \\
k_j(R,S)\le 0
\qquad \text{implies}\qquad
 0\ge k_j(H\cdot R,H\cdot S)\ge k_j(R,S)
\end{align*}

{\sc Remark.} In the case of Riemannian symmetric spaces,
there exists an inverse theorem about characterization of
all angle-compressive maps \cite{Ner-krein}.
 I do not know such statements
 in the $p$-adic case.

\medskip

{\bf 2. Proof}

\smallskip

{\bf 2.1. Norms.}
A {\it norm} on $\Q_p^n$ is a function $v\mapsto\|v\|$
from $\Q_p^n$ to the set
$$\{\dots,p^{-1},1,p,p^2,\dots\}$$
satisfying the conditions

\smallskip

1. $\|v+w\|\le\max(\|v\|,\|w\|)$

2. $\| \lambda v\|=|\lambda| \|v\|$
for $\lambda\in\Q_p$, $v\in \Q^n_p$.

\smallskip

The set of all $v$, such that $\|v\|\le 1$, is a lattice.

Conversely, let $L$ be a lattice. We define the corresponding
norm $\|\cdot\|_L$ by the rule
$$
\|v\|_L:=\Bigl(\min\limits_{k:\,\,p^kv\in L} p^k\Bigr)^{-1}
$$

{\bf 2.2. Minimax characterization of the complex distance.}
The following statement is an exact analogue of
the standard minimax
characterization of the eigenvalues
(see, for instance, \cite{Lid})
and of the stationary
angles (see, for instance \cite {Ner-jordan}.
The proofs also is similar.

\smallskip

{\sc Lemma 1.} {\it  For $L$, $M\in \L_n$}

\begin{align*}
\text{a)}\quad
p^{k_j(R,S)}=
\max\limits_{W:\,\,\dim W=j,\,\,W\subset\Q_p^n}
\Bigl\{
\min\limits_{v\in W,\, v\ne 0} \frac{\|v\|_S}{\|v\|_R}
\Bigr\}
\\
\text{b)} \quad
p^{k_j(R,S)}=
\min\limits_{W:\,\, W\subset \Q_p^n,\,\,\codim W==j-1}
 \Bigl\{
\max\limits_{v\in W,\, v\ne 0} \frac{\|v\|_S}{\|v\|_R}
\Bigr\}
\end{align*}

{\sc Proof.} Let $R$, $S$ have the form (0.2).
Consider the subspace
$Y=\Q_p f_j\oplus\dots\oplus \Q_p f_n$.
We have $\|v\|_R/\|v\|_S\le p^{k_j}$ for each $v\in Y$

 If $\dim W_j=j$, then $Y\cap W\ne 0$, and hence
        $\min_{v\in W}\|v\|_R/\|v\|_S\le p^{k_j}$.

For $W=\Q_p f_1\oplus\dots\oplus \Q_p f_j$ we  have
        $\min_{v\in W}\|v\|_R/\|v\|_S= p^{k_j}$.
\hfill $\square$

\smallskip

{\bf 2.3. Another characterization of the complex distance.}
%
The following statement is obvious.

\smallskip

{\sc Lemma 2.} {\it Let $L$, $M\in\L_n$. Then we have the
following isomorphisms of finite Abelian $p$-groups:}
\begin{align*}
M/(L\cap M)=\bigoplus\limits_{j:\,k_j(L,M)>0} \Z/p^{k_j(L,M)}\Z
\\
L/(L\cap M)=\bigoplus\limits_{j:\,k_j(L,M)<0} \Z/p^{-k_j(L,M)}\Z
\end{align*}

{\bf 2.4. Duality and complex distance.}
Denote by $(\Q_p^n)^\square$ the linear space dual to $\Q_p^n$.
For a lattice $L\in\L_n$ we consider the
{\it dual lattice} $L^\square\subset(\Q_p^n)^\square$
consisting of linear functionals
$f\in(\Q_p^n)^\square$ such that
$$v\in L\,\, \Longrightarrow \,\, f(v)\in\O_p$$

{\sc Lemma 3.}
$
k_j(L,M)=k_j(M^\square,L^\square)
$.

{\bf 2.5. Complex distance and intersections of lattices.}

\smallskip

{\sc Lemma 4.} {\it Let $L$, $M$, $N\in\L_n$. Then
for all $j$, the number $k_j(L\cap M,L\cap N)$ lies between
 $0$ and $k_j(M,N)$.}

\smallskip

This statement is a corollary of Lemma 1
and the following Lemma 5.

\smallskip

{\sc Lemma  5.} {\it For $v\in\Q_p^n$, the number
$\|v\|_{L\cap M}/\|v\|_{L\cap N}$ lies
between 1 and $\|v\|_M/\|v\|_N$.}

\smallskip

{\sc Proof.} Denote
$$
\log_p\|v\|_L=l,\quad
\log_p\|v\|_M=m,\quad
\log_p\|v\|_N=z,\quad
$$
Then
$$
\log_p \|v\|_{L\cap M}=\max (l,m),
\quad
\log_p \|v\|_{L\cap N}=\max (l,z)
$$
It remains to observe that
$\max(l,m)-\max(l,z)$ lies between $0$ and $l-z$.

\smallskip

{\bf 2.6. Complex distance and sums of lattices.}

\smallskip

{\sc Lemma 6.} {\it Let $L$, $M$, $N\in\L_n$.
Then $k_j(L+M,L+N)$ lies between $0$ and $k_j(M,N)$.}

\smallskip

{\sc Proof.} For $P$, $Q\in\L_n$, we have
$$(P+Q)^\square=P^\square\cap Q^\square$$
By the duality consideration, our lemma is equivalent to Lemma 4.

\smallskip

{\bf 2.7. Notation.}
Let $S\subset \Q_p^n\oplus\Q_p^n$ be an element
of the Nazarov semigroup. We define

\smallskip

--- its {\it kernel} $\Ker S=S\cap (\Q_p^n\oplus 0)$,

\smallskip

--- the {\it domain} $\Dom S$ be the image of $S$ under
the projection $\Q_p^n \oplus\Q_p^n\to \Q_p^n\oplus 0$,

\smallskip

--- the {\it image} $\Im S$ be the image of $S$ under the projection
$\Q_p^n \oplus \Q_p^n \to 0\oplus \Q_p^n$,

\smallskip

--- the {\it indefiniteness} $\Indef S$ be the intersection
$S\cap 0\oplus \Q_p^n$

\smallskip

Evidently, all the sets $\Dom S$, $\Ker S$, $\Indef S$, $\Im S$
are lattices in $\Q_p^n$.

\smallskip

{\bf 2.8. Proof of Theorem.}
Let $S\in \Gamma_n$. We have
$$
S\subset \Dom S\oplus  \Im S,\qquad
S\supset \Ker S\oplus \Indef S
.$$
Consider the (finite) factor-group
$$
Y:=\bigl(\Dom S/\Ker S\bigr)\oplus \bigl( \Im S /\Indef S\bigr)
$$
and the image $\sigma$ of $S$ in $Y$.
It can be easily be checked that $\sigma$ is a graph
of an isomorphism
$$
\Dom S/\Ker S\to  \Im S /\Indef S
$$
In particular, these two groups are isomorphic.

Hence, by Lemma 2,
$$k_j(\Dom S, \Ker S)=k_j(\Im S,\Indef S)$$
This implies an existence of a (non canonical)
 element $g_S\in \GL_n(\Q_p)$
such that
$$
g_S\cdot\Dom S=\Im S,\qquad g_S\cdot \Ker S=\Indef S
$$

For any $L\in\L_n$, we have
$$
S\cdot L= g_S \cdot(L\cap \Dom S)+\Indef S
$$
It remains to apply Lemma   4 and 6.

\tt

Math. Phys. Group,

Institute for Theoretical and Experimental Physics,

B. Cheremushkinskaya, 25,  Moscow -- 117259, Russia

\&

University of Vienna

neretin@mccme.ru,


\medskip




\begin{thebibliography}{cc}


\bibitem{Bro}
Brown, K., {\it Buildings,} Springer, 1989.

\bibitem{Ful}
  Fulton, W., {\it Eigenvalues of sums of Hermitian matrices
 (after A. Klyachko)}.
 Asterisque No. 252 (1998), Exp. No. 845, 5, 255--269.


\bibitem{Gar}
Garret, P. {\it Buildings and classical groups,} 1996


\bibitem{Gui}
 Guivarch, Y., Ji, L., Taylor, J. {\it Compactifications
 of symmetric spaces.}
 Birkhauser Boston, Inc., Boston, MA, 1998.

\bibitem{Kou}
Koufany, Kh.,
{\it Contractions of angles in symmetric cones.}
Publ. Res. Inst. Math. Sci. 38 (2002), no. 2, 227--243.



\bibitem{Lid}
Lidskii, V.B. {\it Inequalities for eigenvalues and singular values,}
addendum to Gantmaher F.R. {\it Theory of matrices},  Second (1966), Third
(1976), Forth (1988)
Russian editions.


\bibitem{Mac}
Macdonald, I.G.
{\it Spherical functions on a group of $p$-adic
type,} Publ.Ramanujan Inst., Madras, 1972


\bibitem{Mac2}
Macdonald, I.G.
{\it Symmetric functions and Hall polynomials,}
2-nd ed., Clarendon Press, Oxford, 1995.



\bibitem{Naz}
Nazarov, M.L. {\it The oscillator semigroup over nonarchimedian field.}
J. Funct. Anal., 128 (1995), 384--438

\bibitem{NNO}{Nazarov, V., Yu. Neretin, Yu,, G. Olshanskii, G.,
{\it
Semi-groupes engendr\'es par la repr\'esentation de
Weil du
groupe symplectique de dimension infinie,}
C. R. Acad. Sci. Paris S\'er. I Math. {\bf 309}(1989),
no. 7,
443-446. }


\bibitem{Ner-boson}
{Neretin, Yu., {\it On a semigroup of
operators
in the boson Fock space,} Funct. Anal. Appl. {\bf
24}(1990),
no. 2, 135-144.}





\bibitem{Ner-book}{Neretin Yu., {\it Categories of
symmetries
and infinite-dimensional groups,}
London Math. Soc. Monographs, N. S., {\bf 16},
Oxford University Press, 1996.}


\bibitem{Ner-hinges}
 Neretin, Yu. A.
{\it Hinges and the
Study-Semple-Satake-Furstenberg-De Concini-Procesi-Oshima boundary.}
 Kirillov's seminar on representation theory, 165--230, Amer. Math. Soc.
   Transl. Ser. 2, 181, Amer. Math. Soc., Providence, RI, 1998.

\bibitem{Ner-krein}
  Neretin, Yu. A.
{\it Conformal geometry of symmetric spaces,
 and generalized linear-fractional Kre\u\i n-Shmulian mappings.}
 Mat. Sb. 190 (1999), no. 2, 93--122;
   transl. in Sb. Math. 190 (1999), no. 1-2, 255--283


\bibitem{Ner-jordan}
 Neretin, Yu.A.
 {\it Jordan angles and triangle inequality
  in Grassmannian manifold,} Geometria Dedicata,
86 (2001), 403--432



\bibitem{Ner-tits}
Neretin, Yu.A. {\it Beta-function of Bruhat--Tits buildings
and deformation of $L^2$ on the space of lattices.}
Mat. Sbornik, 2003;
preprint version is available via
{\tt http://xxx.arxiv.org/abs/math/0306079.}


\bibitem{Wei-groups}
Weil, A.
{\it Sur certains groupes d'operateurs unitaires,}
Acta Math., 111 (1964), 143--211


\bibitem{Wei-number}
Weil, A. {\it Basic number theory,} Springer, 1967;

\end{thebibliography}
\end{document}